\documentclass[twocolumn]{autart}

\newtheorem{thm}{Theorem}
\newtheorem{defn}{Definition}
\newtheorem{rem}{Remark}
\newtheorem{exmp}{Example}
\newtheorem{lem}{Lemma}
\newtheorem{prop}{Proposition}
\newtheorem{cor}{Corollary}

\usepackage[sort&compress]{natbib}
\usepackage{amsmath,amssymb}
\usepackage[usenames]{color}
\usepackage{xcolor}
\usepackage{caption,subcaption}
\usepackage{graphicx}
\usepackage{enumerate}
\usepackage{nicefrac}

\usepackage{hyperref}
\hypersetup{
	colorlinks = true,
	allbordercolors = {white},
	allcolors = {blue},
	draft=false
}

\def\ds{\displaystylse}

\def\vec{\mathop{\rm vec}\nolimits}

\DeclareMathOperator*{\argmax}{argmax} 

\def\qed{\hfill$\blacksquare$}

\def\1{\boldsymbol{1}}
\newcommand{\ov}{\overline}
\newcommand{\mc}{\mathcal}
\newcommand{\mb}{\boldsymbol}
\newcommand{\de}{\mathrm{d}}
\newcommand{\0}{\mathbf{0}}

\newcommand{\R}{\mathbb{R}}

\newcommand{\ba}{\begin{array}}
	\newcommand{\ea}{\end{array}}
\renewcommand{\ds}{\displaystyle}
\newcommand{\beq}{\begin{equation}}
	\newcommand{\eeq}{\end{equation}}
\newcommand{\beqn}{\begin{equation*}}
	\newcommand{\eeqn}{\end{equation*}}
\newcommand{\be}{\begin{equation}}
	\newcommand{\ee}{\end{equation}}

\begin{document}
	
	\begin{frontmatter}
		
		\title{Network Behavioral-Feedback SIR Epidemic Model} 
		
		\thanks[footnoteinfo]{Corresponding author: Martina Alutto. \\
			This work was partially supported by the Research Project PRIN 2022 ``Extracting Essential Information and Dynamics from Complex Networks'' (Grant Agreement number 2022MBC2EZ) funded by the Italian Ministry of University and Research.}
		
		\author[KTH]{Martina Alutto \thanksref{footnoteinfo}}\ead{alutto@kth.se},    
		\author[Torino]{Leonardo Cianfanelli}\ead{leonardo.cianfanelli@polito.it},               
		\author[Torino,Lund]{Giacomo Como}\ead{giacomo.como@polito.it},
		\author[Torino]{Fabio Fagnani}\ead{fabio.fagnani@polito.it}  
		
		\address[KTH]{Department of Decision and Control Systems, School of Electrical Engineering and Computer Science, \\KTH Royal Institute of Technology, Stockholm, Sweden.}       
		\address[Torino]{Department of Mathematical Sciences ``G.L.~Lagrange,'' Politecnico di Torino, 10129 Torino, Italy.}  
		\address[Lund]{Department of Automatic Control, Lund University, 22100 Lund, Sweden.}             


		\begin{abstract}   
			We propose a network behavioral-feedback Susceptible-Infected-Recovered (SIR) epidemic model in which the interaction matrix describing the infection rates across subpopulations depends in feedback on the current epidemic state. This model captures both heterogeneities in individuals’ mixing, contact frequency, aptitude to contract and spread the infection, and endogenous behavioral responses such as voluntary social distancing or the adoption of self-protective measures. 
			We first show that the set of equilibrium points coincides with the disease-free states and we analyze their stability.
			Then, we illustrate through several examples how the shape of the set of stable equilibrium points depends on the structure of the interaction matrix, providing insights for the design of effective control strategies. Finally, we analyze the transient dynamic behavior of the epidemics, showing that, in a special case of rank-1 interaction matrices, there always exists an aggregate infection curve that exhibits a unimodal behavior: this result extends the ones on unimodality of the infection curve already known in the literature of epidemic models and paves the way for future control applications.
		\end{abstract}

	\end{frontmatter}
	
	\section{Introduction} 
	Mathematical models play a central role in forecasting the spread of infectious diseases and guiding effective public health interventions. Among these, compartmental models are widely used for the analysis and control of epidemic outbreaks. 
	One of the most studied compartmental models is the Susceptible-Infected-Recovered (SIR) epidemic model, see \cite{Kermack.McKendrick:1927, Hethcote2000TheMO, 
		Diekmann2000MathematicalEO}. In the SIR epidemic model the population is split into three compartments: \textit{susceptible} individuals who have not yet been infected, \textit{infected} individuals who carry the pathogen and can transmit the disease, and \textit{recovered} individuals, who have healed and gained immunity or have otherwise been removed from the system.  
	The classical SIR model has served as a foundation for many recent control strategies, including those developed in response to the COVID-19 pandemic, see, e.g.,  \cite{Cianfanelli.ea:2021, Miclo.ea:2022}.
	
	However, the classical SIR model relies on homogeneity assumptions regarding the pattern of contacts and the transmission rate that are often unrealistic. To address these limitations, network-based epidemic models have been proposed and analyzed  over the last years, see, e.g., \cite{pastor2015epidemic, pare2020modeling, Zino.Cao:2021}. In these models, the nodes of the network represent different subpopulations of indistinguishable individuals that may differ in several aspects, such as geographical area, age, or degree, e.g., in \cite{,hethcote1978immunization, Fall:2007, Nowzari.ea:2016, Mei.ea:2017, Ogura2016StabilityOS,wang2024dynamical}. These subpopulations interact through a network, whose links represent the pattern of contacts. The \textit{interaction matrix} encodes the infection rate across the subpopulations, taking into account factors such as different susceptibility and infectivity levels of the subpopulations, as well as different interaction rates among them.
	
	A large part of the literature on networked epidemic models has focused on the Susceptible-Infected-Susceptible (SIS) epidemic model, where the individuals do not acquire immunity after recovering from the disease. The seminal work by \cite{lajmanovich1976deterministic} established a key bifurcation result: when the largest eigenvalue of the interaction matrix is below a threshold, the unique disease-free equilibrium is stable; otherwise, the system admits an endemic equilibrium. Many extensions of the network SIS model have been then studied, including models with higher-order interactions (\cite{cisnerosvelarde2021multigroup}) and competing viruses (\cite{anderson2023equilibria, liu2019analysis}).
	
	The \textit{network SIR epidemic model} has been the subject of both numerical and theoretical studies, due to its relevance in capturing the spread of epidemics across subpopulations, while accounting for recovery-induced immunity. Network SIR epidemic models have been used to evaluate numerically the effectiveness of targeted interventions, for instance to demonstrate that age-specific lockdown policies outperform uniform lockdown policies in reducing mortality and sustaining economic activity, even in simplified settings with three nodes only and uniform interaction matrix (\cite{acemoglu2021optimal}).  Several works have investigated the theoretical aspects of the network SIR model, see, e.g., \cite{Mei.ea:2017,ellison2020implications,Alutto.ea:2024}. 
	
	A significant part of the literature on the network SIR epidemic model focuses on the case of rank-$1$ interaction matrix, see, e.g., \cite{NOLD1980227,ellison2020implications,Alutto.ea:2024,wang2024dynamical} and \cite[Chapter 3]{hethcote2014gonorrhea}. This case corresponds to assuming that the interaction rate between two subpopulations is proportional to the activity rates of the two subpopulations: as such, it does not account for homophily, that is the tendency of individuals to preferentially interact with others who are similar or belong to the same subpopulation. In particular, \cite{ellison2020implications} studied the special case of symmetric rank-$1$ interaction matrix, that amounts to assuming uniform size, infectivity and susceptibility level across the subpopulations: among other contributions, in this special case \cite{ellison2020implications} provides  a characterization of the set of stable equilibrium points of the model that depends on the heterogeneity of the activity rates of the subpopulations. The analysis of the network SIR model has been extended in \cite{Alutto.ea:2024}, where a number of results are proved, including the identification of $n$ invariants of motion, conditions for the stability of the equilibrium points, and, for the case of rank-$1$ asymmetric interaction matrix with $n$ nodes, a detailed analysis of transient dynamics at the single node level revealing the emergence of multimodal behaviors of the infection curves of single nodes, and the existence of a unimodal weighted aggregate infection curve.

	Besides the assumption of homogeneous population, another key limitation of the classical SIR model is the assumption of a constant transmission rate. This assumption rules out the possibility to model the fact that individuals typically react to the spread of the epidemic and adapt their behavior accordingly (\cite{crosby2003america,lau2005sars}). An approach to overcome this limitation is to consider \textit{behavioral-feedback epidemic models}, incorporating the human endogenous behavior by letting the infection rate depend in feedback on the current epidemic state (\cite{funk2010modelling, verelst2016review}). The first feedback epidemic model was introduced in \cite{CAPASSO197843}, where the authors showed that, if the infection rate is a decreasing function of the fraction of infected individuals, the infection curve remains unimodal as in the classical SIR model. Subsequent studies demonstrated analytically that neglecting the human endogenous behavior may lead to errors in forecasting the trajectory of the epidemic, as the presence of feedback mechanisms can reduce the peak of the infection curve (\cite{Baker2020, gao2024final, franco2020feedback}). Some models extended the analysis by allowing the transmission rate to depend on both the fraction of susceptible and infected individuals (\cite{korobeinikov2006lyapunov}). In particular, \cite{Alutto2021OnSE} proved that the unimodality persists even when the infection rate depends jointly on the fraction of susceptible and infected individuals, if the infection rate is non-decreasing in the fraction of susceptible individuals. 
	
	An alternative approach to incorporate the human behavior is to introduce additional dynamic variables, e.g., to include strategic decision-making such as whether to vaccinate or adopt social distancing based on perceived risk and peer influence (\cite{martins2023epidemic, certorio2022epidemic}), or to integrate opinion dynamics to model how the individuals' beliefs about the epidemic influence adherence to protective measures (\cite{she2022networked, Xu2024,Xuan2020OnAN, paarporn2017networked}), or to model how the infection rate depends on the perceived infection prevalence (\cite{zhou2020actsis, bizyaeva2024active}). The infection rate then depends on these additional variables, whose dynamics is coupled with the epidemic. The literature has mainly focused on the stability of the equilibria and the asymptotic behavior of these models, e.g., to show emergence of oscillations that are not observed in standard epidemic models (\cite{frieswijk2022meanfield, paarporn2023sis, ye2021game, amaral2021epidemiological, kabir2020evolutionary}).

	The analysis of SIR epidemic models that incorporate both network effects and behavioral adaptations remains limited and is primarily empirical (\cite{wang2025, bisin2022spatial, nemati2025modeling}), with the exception of \cite{wang2024dynamical}.
	In that work, the authors introduce a variant of the network SIR model with newborns and deaths, where the incidence and recovery rates depend on feedback from the epidemic state. Individuals are divided into subpopulations based on their degree, and no homophily is considered, so that the resulting interaction matrix is rank-$1$. They further assume that the feedback affects all entries of the matrix uniformly, an assumption we do not make throughout this paper.
	
	The contribution of our paper is three-fold. First, we propose a network SIR model with behavioral-feedback mechanism (NBF-SIR model, in short), which incorporates both network effects and behavioral-feedbacks. These two features were considered only separately in our previous works (\cite{Alutto.ea:2024,Alutto2021OnSE}). In this model, the interaction matrix depends in feedback on the current epidemic state, capturing adaptive contact patterns and heterogeneous behavioral responses across subpopulations, such as voluntary social distancing, behavioral fatigue, or policy-driven restrictions.
	Second, we derive conditions for the stability of the equilibrium points, and explore how different functional forms of the interaction matrix influence the shape of the set of stable equilibrium points. These findings underscore the importance of incorporating behavioral feedbacks in epidemic models and offer insights for the design of effective and targeted control strategies that aim to steer the system towards stable equilibrium points while minimizing the overall fraction of individuals that are infected throughout the outbreak. Note that, in contrast with \cite{wang2024dynamical}, no assumption on the rank of the interaction matrix is made at this level. Third, for a special class of directed rank-$1$ interaction matrices, we establish the existence of a weighted aggregate infection curve that exhibits a unimodal behavior, generalizing the unimodality of the aggregate infection curve established in \cite{Alutto.ea:2024} for the rank-$1$ network SIR model with constant interaction matrix. Since the unimodality of the infection curve is at the basis of the solution of several optimal control problems (\cite{Miclo.ea:2022,CDC2025,alutto2025optimal}), this result may pave the way for control applications, which are left for future research. 
	
	The rest of the paper is organized as follows. In Section \ref{sec:2} we introduce the NBF-SIR model and provide some preliminary results. In Section \ref{sec:3} we study the stability of the equilibrium points and analyze the shape of the set of stable equilibrium points for different interaction matrices. Section \ref{sec:transient} focuses on the transient behavior of the dynamics proving that, for a special class of rank-1 interaction matrices, the dynamics exhibits a unimodal aggregate infection curve. Finally, Section \ref{sec:conclusion} summarizes the work and discusses future research lines.
	\subsection{Notation}
	We briefly gather here some notational conventions adopted throughout the paper. We denote by $\R$ and $\R_{+}$ the sets of real and nonnegative real numbers, respectively. The all-1 vector and the all-0 vector are denoted by $\1$ and $\0$, respectively, whose size may be deduced from the context.
	The spectrum of a matrix $A$ is denoted $\sigma(A)$, and the transpose matrix is denoted $A^T$. For an irreducible matrix $A$ in $\R_+^{n\times n}$, we let $\lambda_{max}(A)$ and $v_{max}(A)$ denote respectively the dominant eigenvalue of $A$ and the corresponding left eigenvector normalized in such a way that $\1^Tv_{max}(A)=1$, {which has positive entries and is unique due to the Perron-Frobenius theorem.} For $x$ in $\R^n$, we let $[x]$ denote the diagonal matrix whose diagonal coincides with $x$. Inequalities between two vectors $x$ and $y$ in $\R^n$ are meant to hold true entry-wise, i.e., $x \le y$ means that $x_i\le y_i$ for every $i$, whereas $x< y$ means that $x_i< y_i$ for every $i$, and $x\lneq y$ means that $x_i\le  y_i$ for every $i$ and $x_j<y_j$ for some $j$. The same is meant to hold true for matrices in $\R^{n\times n}$. The real part of a complex number $x$ is denoted by $\Re(x)$. The closure of a set $\mc S$ is denoted $\bar {\mc S}$.
	
	\section{Model definition}\label{sec:2}
	In this section, we introduce the \emph{network behavioral-feedback SIR epidemic model}.
	We consider a system of $n\ge1$ interacting subpopulations of indistinguishable individuals. We let $x_i(t)$, $y_i(t)$, and $z_i(t)$ denote respectively the fraction of susceptible, infected, and recovered individuals in each subpopulation $i = 1,\cdots,n$ at time $t$, so that 
	\be\label{normalization}x_i(t)+y_i(t)+z_i(t)=1\ee for every subpopulation $i$ and time $t$. 
	We stack these fractions into vectors $x(t)=(x_i(t))$, $y(t)=(y_i(t))$, and $z(t)=(z_i(t))$. 
	Because of the standing validity of relations \eqref{normalization}, we shall denote by $(x,y)$ the state of the system and omit $z$. 
	
	The rate of new infections in each subpopulation due to interactions with infected individuals of all subpopulations is described by the \emph{interaction matrix} $A: [0,1]^{2n} \to \R_+^{n \times n}$, assumed to be a $\mc C^1$ function of the epidemic state $(x,y)$.
	The $(i,j)$-th entry of the interaction matrix $A_{ij}(x, y)$ accounts for the infection rate from subpopulation $j$ to subpopulation $i$, which is modulated by the current epidemic state. This dependence allows to model several behavioral adaptations, e.g., a reduction of interactions with highly infected groups, thereby embedding a feedback mechanism between epidemic prevalence and social interaction patterns. The novelty of this model compared to \cite{Mei.ea:2017} and \cite{Alutto.ea:2024} is that the interaction matrix is not constant.
	Finally, $\gamma>0$ models the recovery rate, which is assumed to be homogeneous across the network.  
	
	The NBF-SIR epidemic model with interaction matrix $A:[0,1]^{2n}\to \R_+^{n\times n}$ and recovery rate $\gamma>0$ is then the autonomous system of ordinary differential equations
	\begin{equation} 
		\label{eq:behavioral-network-SIR}
		\begin{cases}
			\dot{x}_i = - x_i\sum_j A_{ij}(x,y) y_j\\
			\dot{y}_i =   x_i \sum_j A_{ij}(x,y) y_j - \gamma y_i\,, \\
		\end{cases}	\quad i=1,\ldots,n\,,
	\end{equation}
	which, more compactly, reads
	\begin{equation}
		\label{eq:behavioral-network-SIR-compact}
		\dot{x} = -   [x] A(x,y) y\,,\qquad 
		\dot{y} =   [x] A(x,y) y - \gamma y\,.
	\end{equation}
	\begin{rem}
		Note that $A_{ij}(x,y)$ is proportional to the rate of new infections produced in subpopulation $i$ caused by infected individuals of subpopulation $j$, hence $A(x,y)$ does not have to be a row-stochastic matrix. This is consistent with other network SIR models, see \cite{Alutto.ea:2024,Mei.ea:2017}.
	\end{rem}
	\begin{rem}
			Note that the pattern of interactions among the subpopulations can be represented by a weighted directed graph $\mathcal{G}(x,y) = (\mathcal{V}, \mathcal{E}(x,y), A(x,y))$, where $\mc V=\{1,\ldots,n\}$ is the set of nodes each representing a subpopulation, $A(x,y)$ is a state-dependent asymmetric interaction matrix, and $(i,j)$ belongs to the link set $\mc E(x,y)$ if and only if $A_{ij}(x,y)>0$. No assumptions on the symmetry, irreducibility, or rank of $A$ are made at this level.
	\end{rem}
	\begin{exmp}\label{ex1}
		Consider the class of models with interaction matrix
		\be\label{eq:rank1}A_{ij}(x,y) = g_i(x_i)f_j(y_j)\,, \quad i,j = 1,\cdots,n\,,\ee
		with $g_i, f_j:[0,1] \to [0,\infty)$, for every $i,j = 1,\cdots,n$. This model extends the network SIR model with asymmetric rank-1 interaction matrix considered in \cite{Alutto.ea:2024}, by letting the interaction matrix $A$ depend in feedback on the epidemic state. In this model, $g_i(x_i)$ accounts for the susceptibility level and activity rate of subpopulation $i$ and $f_j(y_j)$ accounts for the infectivity level, activity rate and size of subpopulation $j$. Hence, the contact frequency between individuals of subpopulation $i$ and those of subpopulation $j$ is proportional to the activity rates of the two subpopulations, modulated by the epidemic state.
		This class of feedbacks is termed \emph{local} to indicate that $g_i$ depends only on the $i$-th component of $x$ and $f_j$ depends only on the $j$-th component of $y$. 
		Note also that if $f$ and $g$ take only positive values, then the corresponding graph is complete for all $(x,y)$.
		
		A special case of interaction matrix \eqref{eq:rank1} is 
		\be\label{eq:feedback}A_{ij}(x,y) = \frac{g_i(x_i)}{1+ \alpha y_j} \quad \forall i,j = 1,\ldots, n,\ee
		where $g_i : [0,1] \to [0,\infty)$ is an arbitrary non-decreasing function and $\alpha \geq0$. This models the idea that individuals of subpopulation $i$ tend to reduce their activity rate as the fraction of individuals of subpopulation $i$ that are currently infected or have been infected increases, because of an increased risk perception. On the other hand, the monotonicity with respect to $y_j$ may capture how individuals in other subpopulations tend to reduce the interactions with individuals of subpopulation $j$ (or individuals of subpopulation $j$ tend to self-regulate their activity rate) as the fraction of infected individuals in subpopulation $j$ increases. The interaction matrix \eqref{eq:feedback} generalizes previous scalar formulations provided in \cite{franco2020feedback, Baker2020} by incorporating a dependence on the fraction of susceptible individuals and network effects.
	\end{exmp}
	
	The following result gathers some basic properties of the NBF-SIR epidemic model. 
	Let
	$$\mc S:=\left\{(x,y)\in[0,1]^{2n}:x+y\le\1\right\}$$
	denote the space of feasible states. 
	
	\begin{prop}\label{prop:basic}
		Consider the NBF-SIR epidemic model \eqref{eq:behavioral-network-SIR-compact} with $\mc C^1$ interaction matrix $A:[0,1]^{2n}\to \R_+^{n\times n}$ and recovery rate $\gamma >0$. Then:
		\begin{enumerate}
			\item[(i)] for every initial condition $(x(0), y(0))$ in $\mc S$, there exists a unique solution $(x(t), y(t))$, which is $\mc C^1$ and belongs to $\mc S$ for every $t$ in $\R_+$; 
			\item[(ii)] for every $i=1,\ldots,n$, $x_i(t)$ is non-increasing for $t\ge0$, and $x_i(0)>0$ if and only if $x_i(t)>0$ for every $t\ge0$; 
			\item[(iii)] for every $i = 1,\cdots,n$, if $y_i(0) > 0$, then $y_i(t)>0$ for every $t>0$; 
			\item[(iv)] the set of equilibrium points in $\mc S$ is $$\mc S^*:=\left\{(x,\0):\,x\in[0,1]^{n}\right\};$$
			\item[(v)] there exists $x^*$ in $\R_+^n$ such that $\0\le x^*\le x(0)$ and 
			\begin{equation}\label{eq:limit-SIR}
				\lim_{t\to+\infty}(x(t),y(t)) = (x,\0) \in \mc S^*\,.
			\end{equation}
		\end{enumerate}
	\end{prop}
	\begin{pf}
		See Appendix \ref{app}.
	\end{pf}

	\section{Stability of the equilibrium points}\label{sec:3}
	In this section, we analyze the stability properties of the NBF-SIR epidemic model \eqref{eq:behavioral-network-SIR-compact}. We first establish our main result, and then discuss its implications.
	
	\begin{thm}\label{theo:stability}
		Consider the NBF-SIR epidemic model \eqref{eq:behavioral-network-SIR-compact} with $\mc C^1$ interaction matrix $A:[0,1]^{2n}\to \R_+^{n\times n}$ and recovery rate $\gamma>0$. Let $(x,\0)$ in $\mc S^*$ be an equilibrium point.  Then: 
		\begin{enumerate}
			\item[(i)]  if \be\label{lambda>}\lambda_{\max}([x]A(x,\0))>\gamma\,,\ee then $(x,\0)$  is unstable; 
			\item[(ii)]  if \be\label{lambda<}\lambda_{\max}([x]A(x,\0))<\gamma\,,\ee then $(x,\0)$ is stable.
		\end{enumerate}
	\end{thm}
	\begin{pf}
		The Jacobian matrix of \eqref{eq:behavioral-network-SIR-compact} evaluated in an equilibrium point $(x^*,\mb 0)$ in $\mc S^*$ is
		\begin{equation}\label{J}
			J(x^*,\0) = \begin{bmatrix}
				0 &\quad \ds -[x^*]A(x^*,\0)  \\[0.4cm]
				0 &\quad \ds [x^*]A(x^*,\0)- \gamma \mathbf I 
			\end{bmatrix}.
		\end{equation}
		Hence, its spectrum is given by
		$$\sigma(J(x^*,\0))=\sigma([x^*]A(x^*,\0)- \gamma \mathbf I)\cup\{0\}\,.$$
		
		(i) If \eqref{lambda>} holds true, we have that $$\max\{\Re(\lambda) :\, \lambda\in \sigma(J(x^*,\0))\} >0\,.$$ Hence, by the linearization theorem, the equilibrium point $(x^*,\0)$ is unstable. 
		
		(ii) If \eqref{lambda<} holds true, we have that  $$\max\{\Re(\lambda) :\, \lambda\in \sigma(J(x^*,\0))\} = 0\,.$$ In this case, 
		we rely on the center manifold method \cite[Section 3.2]{Guckenheimer1983}. 
		To this end, observe that the center subspace for the equilibrium point $(x^*,\mb 0)$ is 
		$\mathcal{V}_C = \{(x, \0) :\, x \in \R^n\}$. Since by Proposition \ref{prop:basic}(iv) the set of equilibrium points is given by $\mc S^*=\{(x, \0) :\, x \in [0,1]^n\}$, the center manifold coincides with $\mc S^*$ itself. The reduced equation on this manifold is simply given by $\dot{x} = 0$, which implies that $x^*$ is a stable equilibrium point for the reduced system. Therefore, the stability of the equilibrium point $(x^*,\0)$ for the NBF-SIR epidemic model \eqref{eq:behavioral-network-SIR-compact} follows from the center-manifold theorem.  
		\qed \end{pf} 
	
		We now define the set of stable equilibrium points $\mc S^{**}$ of the NBF-SIR model~\eqref{eq:behavioral-network-SIR-compact} as
		$$\mc S^{**} := \left\{ (x^*, \0) \in \mc S^* : \eqref{lambda<}\right\}\,,$$
		whose stability is guaranteed by Theorem~\ref{theo:stability}.
	\begin{rem}
			Note that all stable equilibrium points in 
			$\mc S^*$ are not locally asymptotically stable, since the set of equilibrium points forms a continuum. Moreover, the set $\mc S^{**}$ is not asymptotically stable itself, as every neighborhood of $\mc S^{**}$ contains unstable equilibrium points.
	\end{rem}

	The next result provides sufficient conditions for the monotonicity of stability condition \eqref{lambda<}. 
	\begin{prop}\label{prop:monotone}
		Consider the NBF-SIR epidemic model \eqref{eq:behavioral-network-SIR-compact} with $\mc C^1$ interaction matrix $A:[0,1]^{2n}\to \R_+^{n\times n}$ and recovery rate $\gamma>0$. If
		\begin{align}
			A_{ij}(x,\mb 0)&+x_i\frac{\partial A_{ij}(x,\mb 0)}{\partial x_i}\geq 0\,, \label{eq:cond1} \\
			&\frac{\partial A_{ij}(x,\mb 0)}{\partial x_k} \geq 0 \label{eq:cond2}
		\end{align}
		holds true for every $(x,\mb 0)$ in $\mc S$ and for every $i, j$ and $k$ such that $k\ne i$, then, given a stable equilibrium $(\tilde x,\0)$, every equilibrium $(x,\0)$ with $x \leq \tilde x$ is stable.
	\end{prop}
	\begin{pf}
		From conditions \eqref{eq:cond1}-\eqref{eq:cond2}, it follows that 
		\be \frac{\partial}{\partial x_i} [x]A(x,\0) \geq 0\,, \quad \forall i = 1,\cdots,n\,. \ee
		Therefore, for every $x \leq \tilde x$, we get
		\be  [x]A(x,\0) \leq [\tilde x]A(\tilde x,\0). \ee
		Since for nonnegative square matrices, if $B \leq C$ entry-wise, then $\lambda_{\max}(B) \leq
		\lambda_{\max}(C)$ (\cite{Meyer2000}), we obtain that  
		\be \lambda_{\max}([x]A(x,\0)) \leq \lambda_{\max}([\tilde x]A(\tilde x,\0))\,. \ee
		This implies that if the equilibrium $(\tilde x,\0)$ is stable, every equilibrium $(x,\0)$ with $x \leq \tilde x$ is also stable.
		\qed\end{pf}
	
	Note that conditions \eqref{eq:cond1}-\eqref{eq:cond2} in the scalar case ($n=1$) reduce to requiring that the reproduction number in all equilibrium points is non-decreasing in the fraction of susceptible individuals.
	

	Under the assumption that the subpopulations have equal size (the definition can be adapted if the subpopulations have different size), the set
	$$
	\mc X := \argmax_{\substack{x \in [0,1]^n: \\ (x,\mb 0) \in \ov S^{**}}} \frac 1n \sum_{j = 1}^n x_j
	$$ 
	identifies the points inside the closure of the set of stable equilibrium points that maximize the fraction of susceptible individuals. Note that $\mc X$ belongs to the boundary of $\mc S^{**}$ by construction. Identifying the sets $\mc S^{**}$ and $\mc X$ is crucial to design effective control strategies. Indeed, in optimal control problems with infinite time horizons, the cost of control policies that do not steer the system towards stable equilibrium points diverges in time (\cite{Cianfanelli.ea:2021}). Hence, a sustainable control policy should drive the system toward a stable equilibrium with as large as possible fraction of susceptible individuals, that is, in $\mc S^{**}$ and as close as possible to $\mc X$, to minimize the overall spread of the infection while ensuring long-term stability. The set $\mc X$ was first informally introduced in \cite{ellison2020implications} to discuss the effects of heterogeneity when a planner aims to design a control that drives the system toward stable equilibrium points that maximize the fraction of susceptible individuals in the network SIR model with no feedback and with rank-1 symmetric interaction matrices.
	The next examples illustrate the sets of stable equilibrium points of several NBF-SIR models.
	
	\begin{figure}
		\subfloat[][]{\includegraphics[width=4cm,height=4cm]{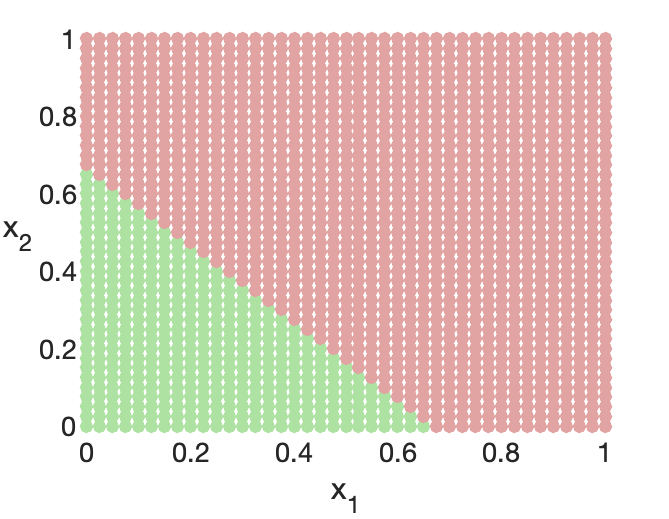}}
		\subfloat[][]{\includegraphics[width=4cm,height=4cm]{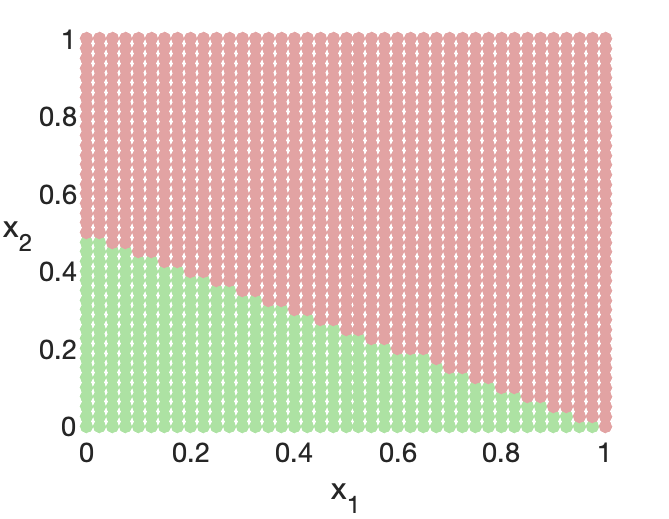}}\\[-0.4cm]
		\subfloat[][]{\includegraphics[width=4cm,height=4cm]{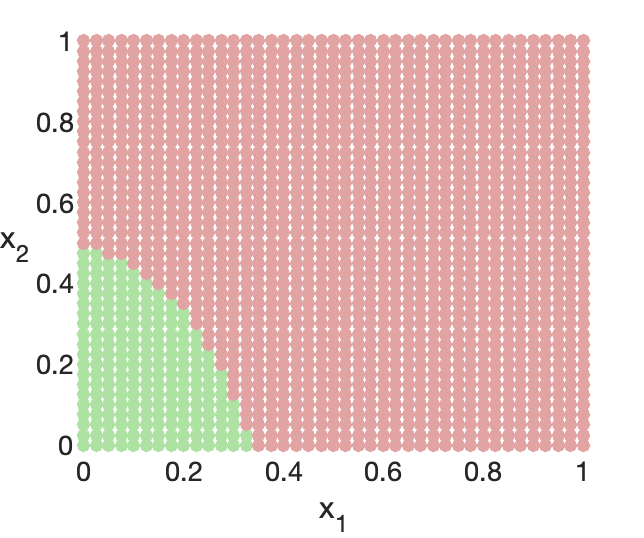}}
		\subfloat[][]{\includegraphics[width=4cm,height=4cm]{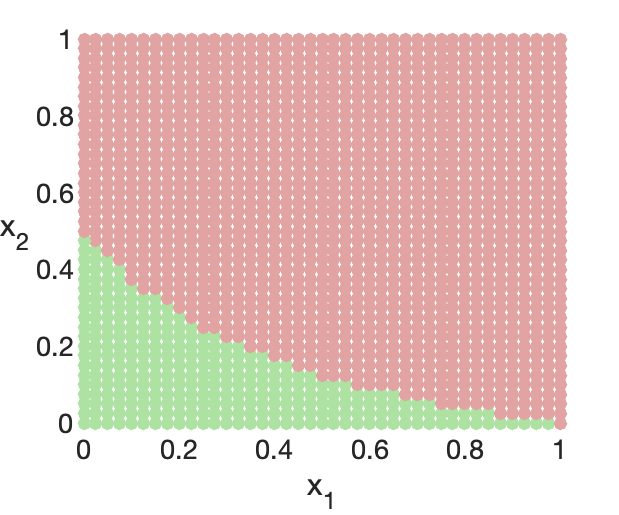}}\\[-0.4cm]
		\subfloat[][]{\includegraphics[width=4cm,height=4cm]{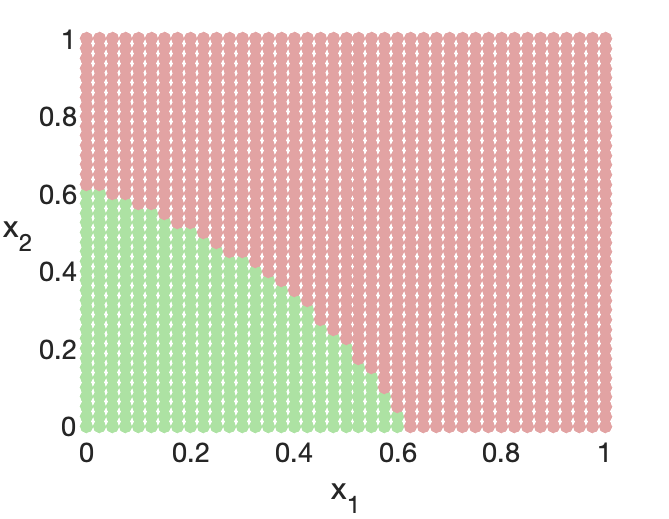}}
		\subfloat[][]{\includegraphics[width=4cm,height=4cm]{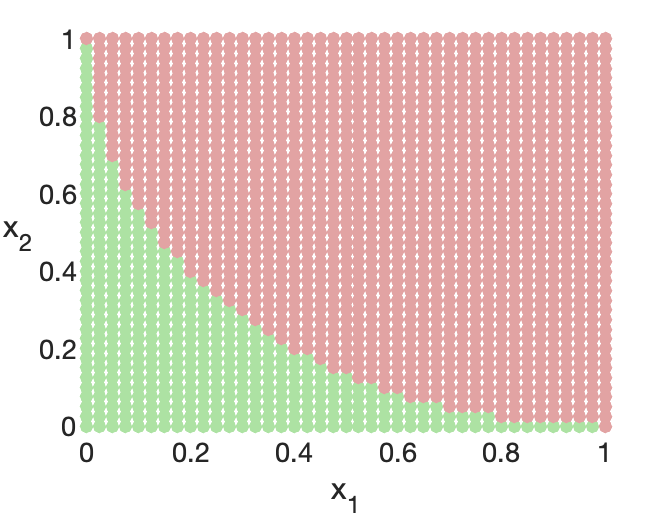}}
		\caption{The set of stable equilibrium points of two-node NBF-SIR epidemic models with different interaction matrices.}
		\label{fig:stability}
	\end{figure}
	
	\begin{exmp}\label{ex:a}
		Consider the NBF-SIR model \eqref{eq:behavioral-network-SIR-compact} with $n=2$ nodes, no feedback, interaction matrix 
		\be\label{eq:A}A = \begin{bmatrix}
			a &\quad b\\
			c &\quad d
		\end{bmatrix}.
		\ee
		and recovery rate $\gamma = 1$. The set of stable equilibrium points is determined by $\lambda_{max}([x]A(x,\mb 0))$, which reads
		$$
		\frac{a x_1 + d x_2}{2} + \sqrt{\bigg(\frac{a x_1 + d x_2}{2}\bigg)^2 - \Delta x_1 x_2}\,,
		$$
		where $\Delta = ad-bc$ is the determinant of the interaction matrix.
		When $A$ is rank-1, i.e., when $\Delta = 0$, the set $\mc S^{**}$ is the intersection between $[0,1]^2$ and a half plane. An example is shown in Figure~\ref{fig:stability}(a) with $A = \frac32 \1 \1^T$. This special case corresponds to an arbitrary partition of a homogeneous population into two subpopulations. 
		In this case, every point on the boundary of $\mc S^{**}$ belongs to $\mc X$, meaning that the planner should drive the system as close as possible to the boundary of $\mc S^{**}$, and all such points are equivalent, which is not surprising, since the model is essentially equivalent to a scalar SIR model. In contrast, Figure \ref{fig:stability}(b) illustrates the set of stable equilibrium points with interaction matrix \eqref{eq:A} and $a=c=1$, $b=d=2$,
			which models two subpopulations with different infectivity levels. Due to this heterogeneity, $\mc X$ is the singleton that contains the maximally asymmetric equilibrium on the boundary of $\mc S^{**}$ where all individuals that have been infected and are now recovered belong to the most infective subpopulation.
			When $A$ is not rank-1, the geometry of $\mc S^{**}$ depends on the sign of $\Delta$. If $\Delta>0$, which is typically associated with homophily (i.e., individuals tend to interact more within their own group), explicit computation shows that $\mc S^{**}$ is a convex region. Conversely, if $\Delta<0$, indicating predominant inter-group interactions, then $[0,1]^{2} \setminus \mc S^{**}$ is convex. These scenarios are illustrated in Figure \ref{fig:stability}(c) for the rank-2 matrix \eqref{eq:A} with $a=3, b=2, c=1, d=2$ and in Figure \ref{fig:stability}(d) for the rank-2 matrix \eqref{eq:A} with $a=1, b=2, c=3, d=2$, respectively. The shape of $\mc S^{**}$ in turn determines the structure of $\mc X$ and the equilibrium points that the planner should target.
		\end{exmp}
		
		
			%
			\begin{exmp}
				Consider the NBF-SIR model \eqref{eq:behavioral-network-SIR-compact} with
				$$A_{ij}(x, y) = \frac{1 + x_i}{1+\alpha y_j}\,, \quad \alpha \geq 0\,, i,j = 1,\cdots,n.$$ 
				This model fits with \eqref{eq:feedback}, as $A(x,y)$ is a rank-1 matrix with local feedback. In this model, the individuals increase their activity rate in response to a higher fraction of susceptible individuals in their subpopulation. Using the fact that $A$ is rank-1, we get that
				$$
				\lambda_{max}([x]A(x,\mb 0)) = \sum_{i = 1}^n x_i(1+x_i)\,,
				$$
				hence the set $\mc S^{**}$ is convex. Figure \ref{fig:stability}(e) illustrates the set of stable equilibrium points with $n=2$, $\alpha = 1.5$ and $\gamma = 1$. The convexity of this region and the symmetry of $A$ imply that $\mc X$ is the singleton that contains the maximally symmetric equilibrium on the boundary of $\mc S^{**}$, therefore the planner should drive the system in $\mc S^{**}$ and as close as possible to such equilibrium.
			\end{exmp}
			
			\begin{exmp}\label{ex:c}
				Consider the NBF-SIR model \eqref{eq:behavioral-network-SIR-compact} with $$A_{ij}(x, y) = \frac{2 - x_i}{1+y_1+y_2},\quad  i,j = 1,\cdots,n.$$ In this case, $A_{ij}(x,y)$ decreases as the fraction of susceptible individuals $x_i$ increases, aligning with the progression of the fraction of individuals of the same population that have been infected. This behavior may represent the effect of increasing pandemic fatigue, where individuals progressively reduce their adherence to restrictive measures as the epidemic evolves and $x_i$ decreases. Moreover, the feedback is non-local, since all entries of the interaction matrix depend on the entire vector $y$, which models situations where the individuals are not aware of the subpopulation that other individuals belong to. In this case,
				$$
				\lambda_{max}([x]A(x,\mb 0)) = \sum_{i = 1}^n x_i(2-x_i)\,.
				$$
				Figure \ref{fig:stability}(f) illustrates the set $\mc S^{**}$ when $n=2$ and $\gamma = 1$.
				The convexity of $[0,1]^2 \setminus \mc S^{**}$ and the symmetry of $A$ imply that $\mc X$ contains the two maximally asymmetric equilibrium points on the boundary of the set $\mc S^{**}$. Hence, the planner should promote lockdown in a single population to drive the system as close as possible to these asymmetric equilibrium points where all individuals that have been infected belong to the same population.
			\end{exmp}
			
			
			%
			
			\begin{rem}
				Note that, since all equilibrium points are disease-free and the stability of an equilibrium does not depend on the derivatives of the interaction matrix, the stability of an equilibrium $(x,\0)$ depends only on $A(x,\0)$. Therefore, the sets of stable equilibrium points illustrated in Figure \ref{fig:stability}(e)-(f) apply to every interaction matrix such that $A_{ij}(x,\0) = 1+x_i$ and $A_{ij}(x,\0) = 2-x_i$, respectively. 
				
			\end{rem}
			
			Examples \ref{ex:a}-\ref{ex:c} show that, depending on the structure of the interaction matrix, the planner should drive the system towards different equilibrium points and design different control strategies accordingly. 
			The analysis applies both to local and non-local feedback mechanisms, to full-rank matrices, and even to asymmetric feedback mechanisms considered in the literature, e.g., a subpopulation (\emph{risk-averters}) reducing their interactions as the infection spreads, similar to those who adhered to strict social distancing measures during the COVID-19 pandemic, and a second subpopulation (\emph{risk-tolerators}) with negative feedback relative to the total fraction of susceptible individuals, see, e.g., \cite{bizyaeva2024active,zhou2020actsis}.

			\section{Transient behavior of the NBF-SIR model}\label{sec:transient}
			In the previous section we have discussed the stability properties of the NBF-SIR model and the implications of network effects and behavioral-feedbacks for the design of control strategies. This section focuses on the transient behavior of the dynamics. The main result of the section is that, under some assumptions on the feedback mechanism, the NBF-SIR model with rank-1 interaction matrix considered in Example \ref{ex1} admits an aggregate infection curve that exhibits a unimodal behavior. 
			Note that the class of NBF-SIR models considered in Example \ref{ex1} can be rewritten as
			\begin{equation}\label{eq:example1}
				\!\!\!\begin{cases}
					\displaystyle	\dot{x}_i = -x_i g_i(x_i) \sum_{j =1}^n f_j(y_j)y_j\\[12pt]
					\displaystyle	\dot{y}_i =   x_i g_i(x_i) \sum_{j =1}^n f_j(y_j) y_j - \gamma y_i
				\end{cases}  \forall i=1,\ldots,n.
			\end{equation}
			For this model, the dominant eigenvalue of $A(x,y)$ is $g(x)^T f(y)$, corresponding to the trace of $A(x,y)$, with associated eigenvector
			$v_{\max}(A(x,y)) = f(y)$. 
			Now, let
			\be \label{bary} \bar{y}(t) = \sum_{j=1}^n f_j(y_j(t))y_j(t) \,.\ee
			The next result states that, under reasonable assumptions, $t\mapsto\bar y(t)$ is unimodal, namely, it exhibits a single peak.

			\begin{thm}\label{theo:bary}
				Consider model \eqref{eq:example1} with initial condition $(x(0),y(0))$ in $\mc S$ such that \be\label{initial-cond} x(0)\gneq\0\,,\qquad y(0)\gneq\0\,.\ee
				Assume that $g(x),f(y)>\0$, that $x_i g_i(x_i)$ is increasing and $f_j(y_j) y_j$ is increasing and concave for every $i,j = 1,\dots,n$ and $(x,y)$ in $\mc S$.
				Then:
				\begin{enumerate} 
					\item[(i)] if $\dot{\bar y}(0) \le 0$, then $\ov y(t)$ is strictly decreasing for $t\ge0$; 
					\item[(ii)] if $\dot{\bar y}(0)>0$, then there exists $\hat t>0$ such that $\ov y(t)$ is strictly increasing on $[0,\hat t]$ and strictly decreasing on $[\hat t,+\infty)$.
				\end{enumerate}
			\end{thm}
			\begin{pf}
				(i)	We prove the result by contradiction. If (i) does not hold true, there must exist $\tau \ge 0$ such that $\dot{\bar y}(\tau)=0$ and $\ddot{\bar y}(\tau) \ge 0$. Direct computation yields
				\be \label{dot-bary} 
				\ba{rcl}\dot{\bar y} & = &  \ds \sum_{j=1}^n \dot y_j [f_j' y_j + f_j]\,,\\[12pt]
				\ddot{\bar y}& = & \ds \sum_{j=1}^n [f_j''y_j + 2 f_j']\dot y_j^2 + \sum_{j=1}^n [f_j'y_j + f_j]\ddot{y}_j
				\ea 
				\ee 
				where, from \eqref{eq:example1},
				\begin{equation}\label{eq:ddoty}
					\ddot y_j = \dot x_j g_j \bar y + x_j g_j' \dot x_j \bar y  +x_j g_j \dot{ \bar y  }- \gamma \dot y_j
				\end{equation}
				for every $j=1,\ldots,n$.
				For simplicity of notation, we omit the dependence on $\tau$ on the right-hand side of the following equation. Note that
				\be\label{eq:ddoty2}
				\ba{rcl}
				\ddot{\bar y}(\tau) & \leq & \ds\sum_{j=1}^n (f_j'y_j + f_j)\ddot{y}_j \\[6pt]
				& \leq & \bar y \ds \sum_{j =1}^n (f_j'y_j + f_j)(g_j' x_j + g_j) \dot x_j - \gamma \dot{\bar y} \\[6pt]
				& = & \bar y \ds \sum_{j =1}^n (f_j'y_j + f_j)(g_j' x_j + g_j) \dot x_j\,,
				\ea
				\ee
				where the first inequality follows from the second equation of \eqref{dot-bary} and from concavity of $f_j(y_j) y_j$, while the second one follows from 
				\eqref{eq:ddoty}, from the first equation of \eqref{dot-bary} and from $\dot{\bar y}(\tau)=0$. Now, note that $x_jg_j(x_j)$ and $f_j(y_j) y_j$ are both increasing for every $j = 1,\dots,n$, hence
				\be\label{eq1}
				\ba{rcl}
				f_j'(y_j(\tau))y_j(\tau) + f_j(y_j(\tau)) & > & 0\,, \\[5pt] 
				g_j'(x_j(\tau)) x_j(\tau) + g_j(x_j(\tau)) & > & 0\,.
				\ea
				\ee 
				Moreover, since $y(0) \gneq \mb 0$ and $f(y)>\0$, Proposition \ref{prop:basic}(iii) implies that 
				\be\label{eq2}
				\bar y(\tau)>0\,.
				\ee 
				Furthermore, 
				\be\label{eq3}
				\dot x_j(\tau) <0\,, \ \ \text{for a } j \in \{1,\cdots,n\}\,.
				\ee
				This follows from $x(0) \gneq \mb 0$ and from Proposition \ref{prop:basic}(ii), which implies that $x_j(\tau)>0$ for at least a $j = 1,\cdots,n$. Hence, $\dot x_j(\tau)<0$ is implied by $\bar y(\tau)>0$, $g(x)>\0$ and \eqref{eq:example1}. Plugging \eqref{eq1}-\eqref{eq3} into \eqref{eq:ddoty2} we get a contradiction and prove the statement. 
				
				(ii) Assume by contradiction that $\dot{\bar y}(t)>0$ for every $t \ge 0$. This would imply the existence of a node $j$ such that $\dot{y}_j(t) > 0$ for all $t \ge 0$, contradicting Proposition \ref{prop:basic}(v). Hence, $\hat t = \min \{ t\geq0 :\dot{\bar y}(t)=0 \}$ is well defined and $\bar y(t)$ is increasing in $[0,\hat t]$. The fact that $\bar y$ is decreasing for $t >\hat t$ follows from point (i).
				\qed\end{pf}
			\begin{rem}
					The unimodality of the infection curve is established by the previous theorem for rank-$1$ interaction matrices. While this assumption may appear restrictive, it is commonly assumed in the literature of network epidemiological models, for instance to study the spread of sexually transmitted diseases (\cite{NOLD1980227}), infectious diseases in populations with heterogeneous activity rates (\cite{ellison2020implications,wang2024dynamical}), or age-structured populations (\cite{acemoglu2021optimal}). In fact, rank-$1$ interaction matrices naturally arise in absence of homophily, whenever the interaction rate between two subpopulations is proportional to the activity rate of the subpopulations. This assumption is reasonable when subpopulations are defined purely by epidemiological attributes (e.g., susceptibility or infectivity levels) rather than by social or spatial proximity. However, it may not be suitable when subpopulations correspond to different geographical areas, where homophily plays a role.
					From a technical viewpoint, the assumption of rank-$1$ interaction matrix simplifies the analysis, since in some cases it allows the infection dynamics to be reduced to a scalar aggregate quantity describing the infection curve on the network, see, e.g., \cite[Chapter 3]{hethcote2014gonorrhea} and (\cite{Alutto.ea:2024}).
			\end{rem}
			\begin{rem}  
				Let
				$h_{ij}(x_i, y_j) = x_i g_i(x_i) f_j(y_j) y_j$
				denote the \emph{force of infection} from node $j$ to node $i$, which represents the rate of new infections in subpopulation $i$ due to interactions with infected individuals of subpopulation $j$.
				The assumptions of Theorem \ref{theo:bary} have a clear epidemiological motivation as they can be equivalently expressed in terms of $h_{ij}(x_i,y_j)$ by saying that $h_{ij}(x_i,y_j)$ is increasing in $x_i$ and $y_j$ and concave in $y_j$. The monotonicity of $h_{ij}(x_i,y_j)$ with respect to $y_j$ reflects the idea that more infected individuals lead to a higher transmission rate, even though the per capita interactions of single infected individuals may decrease due to behavioral adaptations. The concavity captures saturation effects, namely, when more individuals are already infected, the impact of an additional one on the rate of new infections decreases. The monotonicity with respect to $x_i$ reflects the idea that a larger susceptible population increases the potential for new infections. Note that interaction matrix \eqref{eq:feedback} satisfies all assumptions of Theorem \ref{theo:bary}. 
				%
			\end{rem}
			\begin{rem}
				The unimodality of $\bar y = \sum b_j y_j$ was established in \cite[Theorem 1]{Alutto.ea:2024} for the network SIR model with constant rank-1 interaction matrix $$A = ab^T$$
				with vectors $a,b > \mb 0$.
				Since for this constant interaction matrix the force of infection $h_{ij}(x_i,y_j) = a_ix_ib_jy_j$ is increasing in $x_i$ and $y_j$ and concave in $y_j$, and since in this case $f(y) = b$, Theorem \ref{theo:bary} generalizes such a result.
			\end{rem}
			\begin{rem}
				The aggregate variable $\bar{y} = \sum_j f_j(y_j) y_j$ is increasing in every component $y_j$ and may be interpreted as an indicator of the overall infection level in the network. In some special cases, e.g., with interaction matrix \eqref{eq:feedback} with $\alpha = 0$ and when the subpopulations have equal size, $\bar{y} = \sum y_j $ is proportional to the fraction of infected individuals in the entire population.
			\end{rem}
			\medskip

\begin{figure}
				\centering
				\includegraphics[width=5cm,height=4cm]{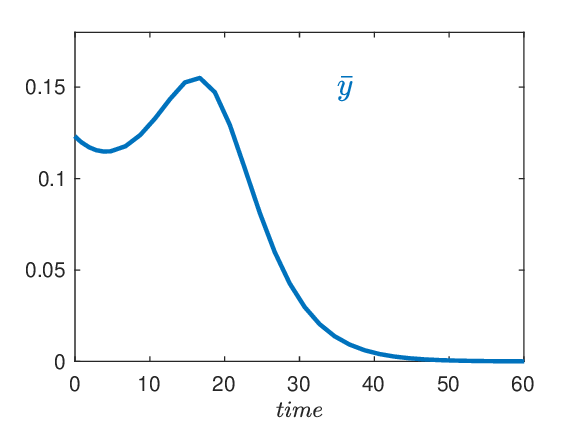}
				\caption{Numerical simulation of the NBF-SIR model of Example \ref{ex:no_uni}.}
			\label{fig:counterexample}
		\end{figure}

			\begin{exmp}\label{ex:no_uni}
				Figure~\ref{fig:counterexample} illustrates a numerical simulation of the NBF-SIR epidemic model \eqref{eq:example1} with $n=5$ subpopulations, interaction matrix $A(x,y) = 0.8(\1 - x)y^T$ and $\gamma = 1$. This interaction matrix is rank-$1$, but $x_ig_i(x_i) = 0.8 x_i (1-x_i)$
				is not increasing in $x_i$ for every $x_i$ in $[0,1]$. As a result, the sufficient conditions established in Theorem~\ref{theo:bary} for the unimodality of the aggregate infection curve are not met. In fact, as illustrated in the figure, the aggregate infection curve exhibits a multimodal behavior. This example highlights that even in the presence of a rank-$1$ interaction matrix, the aggregate infection curve may in general exhibit multimodal behaviors.
			\end{exmp}
			
		\section{Conclusion}\label{sec:conclusion} 
		In this work, we introduced the network behavioral-feedback SIR epidemic model (NBF-SIR), a networked SIR model where the interaction matrix depends in feedback on the epidemic state, to incorporate behavioral adaptations such as voluntary social distancing. We characterized the stability of the equilibrium points, and analyzed how the structure of the interaction matrix affects the set of stable equilibrium points. We then considered the transient dynamics for rank-1 interaction matrices, showing that, under certain assumptions on the feedback mechanism, the dynamics admits a unimodal aggregate infection curve.  
		
		Future research directions include the analysis of the transient dynamics for matrices with higher rank, e.g., to include homophily and heterogeneous interactions, and the design of optimal control strategies for the NBF-SIR model that exploit the properties of the model that have been established in this paper, extending our previous findings regarding optimal control strategies in scalar behavioral-feedback SIR models (\cite{CDC2025,alutto2025optimal}).
		
		\appendix
		\section{Proof of Proposition \ref{prop:basic}}\label{app}
		(i) Since $A$ is a $\mc C'$ function of the state $(x,y)$, and the state space $\mc S$ is compact, existence and uniqueness of the solution are standard  \cite[Chapter~I.3]{hale}. 
		That $\mc S$ is positively invariant simply follows from the fact that for every point belonging to the boundary of $\mc S$, the vector field is either tangent, or points inside the set $\mc S$. 
		
		(ii) Since $\mc S$ is positively invariant, $x(t)$, $y(t)$ and $A(x(t),y(t))$ have non-negative entries for $t\ge0$. Therefore, 
		$$\dot{x}_i(t) = -x_i(t) \sum_{j=1}^n A_{ij}(x,y)y_j(t) \leq 0$$ for all $i = 1, \dots, n$ and $t \geq 0$. By integrating the first equation of \eqref{eq:behavioral-network-SIR}, for every $i=1,\ldots,n$ we get
		$$ x_i(t) = x_i(0) \exp \Big(- \int_{0}^{t} \sum_{j=1}^n A_{ij}(x(\tau), y(\tau)) y_j(\tau) \de \tau \Big),$$
		hence $x_i(t)>0$ for every $t\ge0$ if and only if $x_i(0)>0$. 
		
		(iii) Since $\dot y_i(t) \ge - \gamma y_i(t)$, Gronwall's inequality implies $y_i(t) \ge y_i(0) e^{-\gamma t} > 0$.  
		
		(iv) $(x^*, y^*)$ in $ \mc S$ is an equilibrium point if and only if 
		\begin{equation*}
			\0 = -   [x^*] A(x^*,y^*) y^*\,,\qquad 
			\0 =   [x^*] A(x^*,y^*) y^* - \gamma y^*\,,
		\end{equation*}
		that is if and only if $y^*=\0$.
		
		(v) For $i=1,\ldots,n$, summing up the first two equations of \eqref{eq:behavioral-network-SIR} gives $\dot{x}_i(t) + \dot{y}_i(t) = -\gamma y_i(t)\leq0$, 
		hence $x_i(t)+y_i(t)$ is a non-increasing function.
		Since $x_i(t)$ and $x_i(t) + y_i(t)$ are both non-increasing functions, they admit limits
		$$x_{i}^* := \lim_{t \to +\infty} x_i(t),\;\;\; \xi_{i}^* := \lim_{t \to +\infty} [x_i(t)+y_i(t)]\,.$$
		As a consequence, the following limit 
		\begin{equation*}
			y_{i}^*= \lim_{t \to +\infty} y_i(t)=\xi_{i}^*-x_{i}^*\geq 0
		\end{equation*}
		exists.	Suppose now by contradiction that $y_{i}^* >0$. This means that there exist $T, \epsilon >0$ such that $y_i(t) > \epsilon$, for all $t \geq T$. 
		Then, we get 
		\begin{align*}
			\xi_{i}^*&=\lim_{t \to +\infty} [x_i(t) + y_i(t) ]\\
			&= x_i(0) + y_i(0) - \gamma \lim_{t \to +\infty} \int_{0}^{t} y_i(\tau) d\tau \\
			&< x_i(0) + y_i(0) - \gamma \lim_{t \to +\infty} \int_{T}^{t} y_i(\tau) d\tau \\[6pt]
			&< x_i(0) + y_i(0) - \gamma \epsilon \lim_{t \to +\infty} (t- T) = - \infty
		\end{align*}
		As this is a contradiction, it must be $y_{i}^* =0$ and this is valid for every $i = 1,\dots,n$. 
		Therefore, for every $(x(0), y(0))$ in $\mc S$, the solution $(x(t), y(t))$ converges to an equilibrium point $(x^*, \0)$. Moreover, from point (ii), we have that $\0 \leq x^* \leq x(0)$.
		\qed
				
				\bibliographystyle{agsm}
				
				\bibliography{bib}
				
			\end{document}